\documentclass[12pt]{article}

\usepackage{amsmath,amssymb,amsthm}
\usepackage{graphicx}

\newtheorem{theorem}{Theorem}
\newtheorem{claim}{Claim}

\newcommand{\li}{{\mathrm{lim\,inf}}}
\newcommand{\ls}{{\mathrm{lim\,sup}}}

\title{Maximum sparse induced subgraphs of the binomial random graph with given number of edges}
\author{Dmitry Kamaldinov\footnote{Moscow Institute of Physics and Technology (National Research University), 9 Institutskiy per., Dolgoprodny, Moscow Region, 141701, Russian Federation}, Arkadiy Skorkin\footnote{Adyghe State University, ul. Pervomayskaya, 208, Maykop, Republic of Adygea, 385000, Russian Federation}, Maksim Zhukovskii\footnote{Moscow Institute of Physics and Technology (National Research University), laboratory of advanced combinatorics and network applications, 9 Institutskiy per., Dolgoprodny, Moscow Region, 141701, Russian Federation; Adyghe State University, Caucasus mathematical center, ul. Pervomayskaya, 208, Maykop, Republic of Adygea, 385000, Russian Federation; The Russian Presidential Academy of National Economy and Public Administration, Prospect Vernadskogo, 84, bldg 2, Moscow, 119571, Russian Federation.}}

\date{}

\begin{document}
\maketitle

\begin{center}
{\bf Abstract}
\end{center}

We prove that a.a.s. the maximum size of an induced subtree of the binomial random graph $G(n,p)$ is concentrated in 2 consecutive points. We also prove that, given a non-negative integer-valued function $t(k)<\varepsilon k^2$, under a certain smoothness condition on this function, a.a.s. the maximum size $k$ of an induced subgraph with exactly $t(k)$ edges of $G(n,p)$ is concentrated in 2 consecutive points as well.

\vspace{1cm}

\section{Introduction}

In~\cite{BE_independent}--\cite{M_independent}, it was proven that a.a.s. (with asymptotical probability 1) the maximum size of an independent set (the independence number) of the binomial random graph $G(n,p=\mathrm{const})$ (see, e.g., \cite{Bollobas}--\cite{ZhukRodi}) equals either $f_0(n)$ or $f_0(n)+1$, where
$$
 f_0(n)=\left\lfloor2\log_{1/(1-p)}n-2\log_{1/(1-p)}\log_{1/(1-p)}n+2\log_{1/(1-p)}\frac{e}{2}+0.9\right\rfloor
$$
(the above remains true after substituting $1-\varepsilon$ instead of 0.9 for an arbitrary positive $\varepsilon<0.5$). Is the 2-point concentration also true, say, for the maximum induced path, maximum induced cycle or the maximum induced tree? Is it true for the maximum subgraph with certain restrictions on the number of edges?\\

In~\cite{Many_Induced_Graphs}, 2-point concentration results were obtained for the simple path and for the simple cycle.

\begin{theorem}[Dutta K., Subramanian C.R., 2018]
Let $f_*(n)=\left\lfloor2\log_{1/(1-p)}np+2.9\right\rfloor$. Then a.a.s. both the maximum size of an induced path and the maximum size of an induced cycle in $G(n,p)$ belongs to $\{f_*(n),f_*(n)+1\}$.
\label{many_induced}
\end{theorem}

In the paper, they also ask about a 2-point concentration result for trees but failed in answering this question. We answer this question below.

\begin{theorem}
There exists an $\varepsilon>0$ such that a.a.s. the maximum size of an induced tree in $G(n,p)$ belongs to $\{f_{\varepsilon}(n),f_{\varepsilon}(n)+1\}$, where
$$
 f_{\varepsilon}(n)=\lfloor2\log_{1/(1-p)}(enp)+2+\varepsilon\rfloor.
$$ 
\label{th_trees}
\end{theorem}

The proof is given in Section~\ref{trees_proof}.\\

In \cite{Subgraphs_bounded_deg}, subgraphs with bounded from above numbers of edges were considered.  The main result of \cite{Subgraphs_bounded_deg} states, in particular, the following. 

\begin{theorem}[Fountoulakis N., Kang R.J., McDiarmid C., 2014]
Let $t(k)$ be a sequence of non-negative integers such that $t=o(\sqrt{k^3/\ln k})$. Then, a.a.s. the maximum size $k$ of an induced subgraph in $G(n,p)$ with at most $t(k)$ edges belongs to the set of 2 consecutive numbers $\{f_t(n),f_t(n)+1\}$, where
$$
 f_t(n)=\left\lfloor2\log_{1/(1-p)}n+(t-2)\log_{1/(1-p)}\log_{1/(1-p)}n-t\log_{1/(1-p)} t	+t\log_{1/(1-p)}(2bpe)+2\log_{1/(1-p)}\frac{e}{2}+0.9\right\rfloor.
$$
\end{theorem}

Moreover, for subgraphs with the number of edges equal to a given function $t(k)=p{k\choose 2}+O(k)$ of the number of its vertices $k$, in~\cite{BaloghZhuk}, the 2-point concentration was disproved.

\begin{theorem}[Balogh J., Zhukovskii M., 2019]
Let $t(k)=p{k\choose 2}+O(k)$ be a sequence of positive integers. Let $X_n$ be the maximum size $k$ of an induced subgraph in $G(n,p)$ with $t(k)$ edges. There exists $\mu>0$ such that, for $c>\mu$ and $C>2c+\mu$, we have
$$
0<\li_{n\to\infty}{\sf P}\left(n-C\sqrt{\frac{n}{\ln n}}<X_n<n-c\sqrt{\frac{n}{\ln n}}\right)\leq
$$
$$
 \ls_{n\to\infty}{\sf P}\left(n-C\sqrt{\frac{n}{\ln n}}<X_n<n-c\sqrt{\frac{n}{\ln n}}\right)<1.
$$
Moreover, let, for any sequence $m_k=O(\sqrt{k/\ln k})$ of non-negative integers, 
$$
\left|\left(t(k)-{k\choose 2}p\right)-\left(t(k-m_k)-{k-m_k\choose 2}p\right)\right|=o(k).
$$
Then, for every $\varepsilon>0$, there exist $c,C$ such that
$$
 \li_{n\to\infty}{\sf P}\left(n-C\sqrt{\frac{n}{\ln n}}<X_n<n-c\sqrt{\frac{n}{\ln n}}\right)>1-\varepsilon.
$$
\end{theorem}

It is natural to ask about the 2-point concentration result for smaller $k$. In Section~\ref{sparse_proof}, we prove the following.

\begin{theorem}
Let $R>0$. There exists an $\varepsilon>0$ such that, for every sequence of non-negative integers $t(k)$ with the following properties
\begin{itemize}
\item $|t(k+1)/t(k)-1|\leq\frac{R}{k}$ for all large enough $k$,
\item $t(k)<\varepsilon k^2$  for all large enough $k$,
\end{itemize}
there exists $f(n)$ such that $|f(n)-2\log_{1/(1-p)}n|\leq (3\varepsilon\ln\frac{1}{\varepsilon})\ln n$ and a.a.s. the maximum size $k$ of a set with $t(k)$ edges belongs to $\{f(n),f(n)+1\}$.
\label{th_edges}
\end{theorem}

In other words, a threshold on the number of edges for the 2-point concentration is $\Theta(k^2)$.\\

{\it Remark 1}. The first smoothness condition in Theorem~\ref{th_edges} can not be removed because, clearly, one may consider two different sequences (say, $t_1(k)=1$, $t_2(k)=\lfloor\varepsilon k^2\rfloor -1$), and combine them in the sequence $t(k)=t_1(k)$ for odd $k$ and $t(k)=t_2(k)$ for even $k$. For such a sequence the 2-point concentrations result fails, but it is true both for $t_1(k)$ and $t_2(k)$. Notice that, say, $t(k)=Ck^{a}(1+O(1/k))$, $a>0$, satisfies the condition.

It is also worth mentioning that $f(n)=2\log_{1/(1-p)}n(1+o(1))$, if $t=o(k^2)$. \\

{\it Remark 2}. For several other random subgraph models, the independence number has another asymptotical behaviour~(see, e.g.,~\cite{DerevyankoKiselev,RaiStability}). It would be of interest to study the maximum size $k$ of an induced subgraph with $t(k)$ edges for these models as well.

\section{Maximum induced trees}
\label{trees_proof}

Let $X_k$ be the number of induced subtrees in $G(n,p)$ of size $k$. Clearly,
$$
 {\sf E}X_k={n\choose k}(1-p)^{{k\choose 2}-k+1}p^{k-1}k^{k-2}.
$$
For $k=O(\ln n)$, we get ${\sf E}X_k\sim e^{k\ln n-\frac{5}{2}\ln k+k-{k\choose 2}\ln[1/(1-p)]+(k-1)\ln[p/(1-p)]-\frac{1}{2}\ln(2\pi)}=:\gamma(k)$.

Since $\frac{\partial\ln\gamma(k)}{\partial k}=\ln n-\frac{5}{2k}+1-k\ln[1/(1-p)]+\ln[p/(1-p)^{3/2}]<0$ for $k>\frac{3}{2}\frac{\ln n}{\ln[1/(1-p)]}$ and large enough $n$, there exists 
$$
\hat k(n)=2\frac{\ln n}{\ln[1/(1-p)]}+O(1)
$$ 
such that $\gamma(k)=1$ and, therefore, for $\varepsilon>0$,
\begin{equation}
 {\sf E}X_{\lceil\hat k+\varepsilon\rceil}\leq \gamma(\hat k+\varepsilon)(1+o(1))= e^{\varepsilon\ln n-\varepsilon\hat k\ln(1/(1-p))+O(1)}=e^{-\varepsilon\ln n+O(1)}\to 0,
\label{expectation_upper}
\end{equation}
\begin{equation}
 {\sf E}X_{\lfloor\hat k-1+\varepsilon\rfloor}\geq\gamma(\hat k-1+\varepsilon)(1+o(1))= e^{-(1-\varepsilon)\ln n+(1-\varepsilon)\hat k\ln(1/(1-p))+O(1)}=e^{(1-\varepsilon)\ln n+O(1)}\to\infty.
\label{expectation_lower}
\end{equation}
By Markov inequality, from (\ref{expectation_upper}), we get ${\sf P}(t(G(n,p))<\lceil\hat k+\varepsilon\rceil)\to 1$.\\

It remains to prove that ${\sf P}(t(G(n,p))\geq\lfloor\hat k-1+\varepsilon\rfloor)\to 1$. 

For this, set $k=\lfloor\hat k-1+\varepsilon\rfloor$. In the usual way, we get the following bound for the second factorial moment:
$$
 {\sf E}X_k(X_k-1)-{\sf E}X_k^2\leq\sum_{\ell=2}^{k-1}F_{\ell},
$$
$$
F_{\ell}= 
 {n\choose k}k^{k-2}{k\choose\ell}{n-k\choose k-\ell}\max_{r\in\{0,\ldots,\ell-1\}}p^{2(k-1)-r}(1-p)^{2{k\choose 2}-{\ell\choose 2}-2(k-1)+r}f(k,\ell,r),
$$
where $f(k,\ell,r)$ is an upper bound (we will define its precise value below) for the number of trees on a set of vertices $\{a_1,\ldots,a_k\}$ such that the set of edges of these trees between the vertices of $\{a_1,\ldots,a_{\ell}\}$ is fixed and has cardinality $r$. \\

For $\ell\leq 2\frac{\ln n}{\ln[1/(1-p)]}-6\frac{\ln\ln n}{\ln[1/(1-p)]}$, we will use the trivial bound $f(k,\ell,r)=k^{k-2}$. For such $\ell$,
$$
 \frac{F_{\ell}}{({\sf E}X_k)^2}=\frac{{k\choose\ell}{n-k\choose k-\ell}(1-p)^{-{\ell\choose 2}}\max_{r\in\{0,\ldots,\ell-1\}}((1-p)/p)^{r}}{{n\choose k}}.
$$
If $1-p\leq p$, then
\begin{equation}
 \frac{F_{\ell}}{({\sf E}X_k)^2}\leq\frac{{k\choose\ell}{n-k\choose k-\ell}(1-p)^{-{\ell\choose 2}}}{{n\choose k}}\leq\sqrt{2\pi k}\left(\frac{k^2e }{\ell n}(1-p)^{-\frac{\ell}{2}}\right)^{\ell}\leq (k^{-1+o(1)})^{\ell}.
\label{raz}
\end{equation}
If $1-p>p$, then
\begin{equation}
 \frac{F_{\ell}}{({\sf E}X_k)^2}\leq\frac{{k\choose\ell}{n-k\choose k-\ell}(1-p)^{-{\ell\choose 2}}((1-p)/p)^{\ell}}{{n\choose k}}\leq\sqrt{2\pi k}\left(\frac{k^2e(1-p)}{p\ell n}(1-p)^{-\frac{\ell}{2}}\right)^{\ell}\leq (k^{-1+o(1)})^{\ell}.
\label{dva}
\end{equation}

Let us switch to $\ell>2\frac{\ln n}{\ln[1/(1-p)]}-6\frac{\ln\ln n}{\ln[1/(1-p)]}$. 

In order to define $f(k,\ell,r)$, consider a tree $T$ with a vertex set $A$ of cardinality $k$. Let $B\cap A=\Upsilon$, $|\Upsilon|=\ell$ and $|B|=k$ as well. Assume that $T$ has exactly $r$ edges in $\Upsilon$.  Let us estimate from above the number of trees on the vertex set $B$ such that, in $\Upsilon$, they induce the same set of edges. Set $B\setminus A:=\{v_1,\ldots,v_{k-\ell}\}$.

Let $\mathcal{H}$ be the set of connected components of $A|_{\Upsilon}$. Clearly, $|\mathcal{H}|=\ell-r$.

There are $(k-\ell)^{\ell-r}$ decompositions of $\mathcal{H}$ into $k-\ell$ parts. Let $\mathcal{H}_1\sqcup\ldots\sqcup\mathcal{H}_{k-\ell}$ be such a decomposition. 

In every component of $\mathcal{H}$, choose a vertex (there are at most $(\frac{\ell}{\ell-r})^{\ell-r}$ ways of doing that for $\ell-r\leq\frac{\ell}{e}$, at most $3^{2r-\ell}2^{2\ell-3r}$ ways for $\frac{\ell}{2}\geq\ell-r>\frac{\ell}{e}$ and at most $2^{r}$ ways for $\ell-r>\ell/2$). Then, for every $i\in\{1,\ldots,k-\ell\}$, and every component of $\mathcal{H}_i$, we draw an edge between the chosen vertex and $v_i$.

Consider the set of $k-\ell$ vertices $U_1,\ldots,U_{k-\ell}$, where $U_i$ is the union of the set of all vertices of all components of $\mathcal{H}_i$ and $\{v_i\}$. There are $(k-\ell)^{k-\ell-2}$  trees on this set. Let $\mathcal{F}$ be such a tree. For every edge $\{U_i,U_j\}$ of this tree, draw, in $B$, an edge $e_{ij}$ between a vertex $u_i\in U_i$ and a vertex $u_j\in U_j$ such that either $u_i$ or $u_j$ is in $\Upsilon$. Clearly,  there are at most $\ell+1$ ways of choosing $e_{ij}$. Clearly, the final graph on the set $B$ is a tree, and every tree can be constructed using the above procedure.

Putting all together, there are at most
$$
f(k,\ell,r):=\left(\frac{\ell}{\ell-r}\right)^{\ell-r}(k-\ell)^{k-r-2}(\ell+1)^{k-\ell-1}
$$ 
trees for $r\geq\ell(1-1/e)$, at most  
$$
f(k,\ell,r):=3^{2r-\ell}2^{2\ell-3r}(k-\ell)^{k-r-2}(\ell+1)^{k-\ell-1}
$$ 
trees for $\ell/2\leq r<\ell(1-1/e)$, and at most
$$
f(k,\ell,r):=2^{r}(k-\ell)^{k-r-2}(\ell+1)^{k-\ell-1}
$$ 
trees for $r<\ell/2$.\\

Set $g(k,\ell,r)=f(k,\ell,r)((1-p)/p)^r$. Then
$$
 \frac{\partial}{\partial r}\ln g(k,\ell,r)=
$$
$$ 
 \ln\left(\frac{e(\ell-r)(1-p)}{\ell(k-\ell)p}\right)I[r\geq\ell(1-1/e)]+\ln\left(\frac{9(1-p)}{8(k-\ell)p}\right) I[\ell/2\leq r<\ell(1-1/e)]+\ln\left(\frac{2(1-p)}{(k-\ell)p}\right)I[r<\ell/2].
$$
If $\ell\leq k-\frac{2(1-p)}{p}$, then $g(k,\ell,r)$ dicreases with $r$. Therefore, its maximum equals $g(k,\ell,0)=(k-\ell)^{k-2}(\ell+1)^{k-\ell-1}$. Then,
$$
 \frac{F_{\ell}}{({\sf E}X_k)^2}\leq\frac{{k\choose\ell}{n-k\choose k-\ell}(1-p)^{-{\ell\choose 2}}(k-\ell)^{k-2}(\ell+1)^{k-\ell-1}}{{n\choose k}k^{k-2}}=:\tilde F_{\ell}.
$$
Consider the fraction
\begin{equation}
\frac{\tilde F_{\ell+1}}{\tilde F_{\ell}}=\frac{(k-\ell)^2}{(\ell+1)(n-2k+\ell+1)}(1-p)^{-\ell}\left(\frac{k-\ell-1}{k-\ell}\right)^{k-2}\left(\frac{\ell+2}{\ell+1}\right)^{k-\ell-2}\frac{1}{\ell+1}.
\label{F_fraction}
\end{equation}
For $n$ large enough,
$$
\frac{\partial}{\partial\ell}\ln\left[\frac{\tilde F_{\ell+1}}{\tilde F_{\ell}}\right]=
$$
$$
-\frac{2}{k-\ell}-\frac{2}{\ell+1}-\frac{1}{n-2k+\ell+1}+\ln\left[\frac{1}{1-p}\right]-\frac{k-2}{(k-\ell)(k-\ell-1)}-\ln\frac{\ell+2}{\ell+1}-\frac{k-\ell-2}{(\ell+1)(\ell+2)}<0,
$$
since $\ell>2\frac{\ln n}{\ln[1/(1-p)]}-6\frac{\ln\ln n}{\ln[1/(1-p)]}$ and $k=2\frac{\ln n}{\ln[1/(1-p)]}+O(1)$. Therefore, $\frac{\tilde F_{\ell+1}}{\tilde F_{\ell}}$ dicreases with $\ell$ in the range.

Let $\ell=k-\max\left\{2,\left\lceil\frac{8}{\ln[1/(1-p)]}\right\rceil\right\}$. Then, for $n$ large enough, from~(\ref{F_fraction}), we get
$$
 \frac{\tilde F_{\ell+1}}{\tilde F_{\ell}}>\frac{(1-p)^{-\ell}}{\ell^2 n}e^{-2k/(k-\ell)}\geq e^{\frac{1}{2}\ln n+O(\ln\ln n)}.
$$
Therefore, $\tilde F_{\ell}$ increases in
$\left(2\frac{\ln n}{\ln[1/(1-p)]}-6\frac{\ln\ln n}{\ln[1/(1-p)]},k-\max\left\{2,\frac{8}{\ln[1/(1-p)]}\right\}\right]$.

Furthermore, for $\ell=k-c$, $c\in\mathbb{N}$, we get
$$
 \tilde F_{\ell}\leq \frac{{k\choose c}{n-k\choose c}(1-p)^{ck-c/2-c^2/2}(p/(1-p))^{k-1}c^{k-2}k^c}{{n\choose k}k^{k-2}(1-p)^{{k\choose 2}}(p/(1-p))^{k-1}}=
$$
$$ 
 \frac{\exp(c\ln n-c k\ln(1/(1-p))+k\ln[p/(1-p)]+k\ln c+o(k))}{{\sf E}X_k}=
$$
$$ 
 \frac{\exp[-k(c\ln(1/(1-p))/2-\ln(p/(1-p))-\ln c)(1+o(1))]}{{\sf E}X_k}.
$$
The function $c\ln(1/(1-p))/2-\ln c$ approaches infinity as $c\to\infty$. Therefore, for $c$ large enough, $c\ln(1/(1-p))/2-\ln(p/(1-p))-\ln c>0$.

Let us conclude the above arguments: there exists $c\in\mathbb{N}$ and $a>0$ such that, for all $\ell\in\left(2\frac{\ln n}{\ln[1/(1-p)]}-6\frac{\ln\ln n}{\ln[1/(1-p)]},k-c\right]$, and $n$ large enough,
\begin{equation}
 \frac{F_{\ell}}{({\sf E}X_k)^2}\leq \frac{e^{-ak}}{{\sf E}X_k}.
\label{medium}
\end{equation}
By the Chebyshev inequality, from (\ref{expectation_lower}), (\ref{raz}), (\ref{dva}) and (\ref{medium}) we get
$$
 {\sf P}(X_k=0)\leq\frac{{\sf D}X_k}{({\sf E}X_k)^2}\leq\frac{\sum_{\ell=2}^{k-1}F_{\ell}+{\sf E}X_k}{({\sf E}X_k)^2}=\sum_{\ell=k-c}^{k-1}\frac{F_{\ell}}{({\sf E}X_k)^2}+o(1).
$$

It remains to prove that, for every $c\in\mathbb{N}$, 
$$
\max_{\ell\in\{k-c,\ldots,k-1\}}\frac{F_{\ell}}{({\sf E}X_k)^2}=o(1).
$$
Let $c\in\mathbb{N}$, $\ell=k-c$. Then, for some constant $A>0$,
$$
 \frac{F_{\ell}}{({\sf E}X_k)^2}\leq \frac{A{k\choose c}{n-k\choose c}(1-p)^{ck}k^c\max_{r\in\{0,1,\ldots,k-c-1\}} f_0(k,r) c^{k-r}(p/(1-p))^{k-r}}{{n\choose k}k^{k-2}(1-p)^{{k\choose 2}}(p/(1-p))^{k-1}},
$$
where 
$$
f_0(k,r)=\left(\frac{\ell}{\ell-r}\right)^{k-r}I(r\geq\ell(1-1/e))+
(4/3)^k (9/8)^rI(\ell/2\leq r<\ell(1-1/e))+2^rI(r<\ell/2).
$$
Therefore,
$$
 \frac{F_{\ell}}{({\sf E}X_k)^2}\leq \frac{\exp\left[-c\ln n+\max\limits_{r\in\{0,\ldots,k-c-1\}}\ln f_1(k,r)+O(\ln\ln n)\right]}{{\sf E}X_k},\quad
 f_1(k,r)=f_0(k,r)c^{k-r}(p/(1-p))^{k-r}.
$$
If $r\geq\ell(1-1/e)$, then 
$$
\frac{\partial\ln f_1(k,r)}{\partial r}=-\ln\frac{pc}{(1-p)}-\ln\frac{\ell}{\ell-r}+1+\frac{c}{\ell-r}.
$$
Therefore, $f_1$ decreases with $r$ if $\frac{pc}{1-p}>1$, and first increases and then decreases if $\frac{pc}{1-p}\leq 1$. If $\ell/2\leq r<\ell(1-1/e)$, then $\frac{\partial\ln f_1(k,r)}{\partial r}=\ln\frac{9(1-p)}{8pc}$. Therefore, $f_1$ increases with $r$ if and only if $\frac{pc}{1-p}<\frac{9}{8}$. Finally, if $r<\ell/2$, $\frac{\partial\ln f_1(k,r)}{\partial r}=\ln\frac{2(1-p)}{pc}$. Therefore, $f_1$ increases with $r$ if and only if $\frac{pc}{1-p}<2$.

Notice that $f_1(k,r)$ has two discontinuities: $r_1=\frac{\ell}{2}$, $r_2=\ell-\frac{\ell}{e}$. It is easy to check that $f_1(k,r_1)>f_1(k,r_1-0)$ and $f_1(k,r_2)>f_2(k,r_2-0)$. Moreover, if $\frac{pc}{1-p}\geq 2$, then $f_1(k,0)>f_1(k,r_1)>f_1(k,r_2)$.

Summing up, 

\begin{itemize}

\item if $\frac{pc}{1-p}\leq 1$, then $f_1$ achieves its maximum in $r\sim \ell(1-\frac{pc}{e(1-p)})$;
\item if $1<\frac{pc}{1-p}\leq \frac{9}{8}$, then $f_1$ achieves its maximum in $r=\ell(1-1/e)$;
\item if $\frac{9}{8}<\frac{pc}{1-p}\leq 2$, then $f_1$ achieves its maximum either in $r=\ell/2$, or in $r=\ell(1-1/e)$;
\item if $\frac{pc}{1-p}>2$, then $f_1$ achieves its maximum in $r=0$.

\end{itemize}

Below, we consider all these four situations separately.

\begin{enumerate}

\item Let $\frac{pc}{1-p}\leq 1$. If $c=1$, then $p\leq\frac{1}{2}$. Therefore,
$$
 \frac{F_{\ell}}{({\sf E}X_k)^2}\leq \frac{\exp\left[-\ln n\left(1-\frac{2p}{e(1-p)\ln[1/(1-p)]}+o(1)\right)\right]}{{\sf E}X_k}\leq e^{-(2-2/[e\ln 2]-\varepsilon+o(1))\ln n}.
$$
If $c\geq 2$, then $p\leq\frac{1}{3}$. Therefore, 
$$
 \frac{F_{\ell}}{({\sf E}X_k)^2}\leq \frac{\exp\left[-c\ln n\left(1-\frac{2p}{e(1-p)\ln[1/(1-p)]}+o(1)\right)\right]}{{\sf E}X_k}\leq \frac{\exp\left[-c\ln n\left(1-\frac{1}{e\ln(3/2)}+o(1)\right)\right]}{{\sf E}X_k}.
$$

\item Let $1<\frac{pc}{1-p}<2$ and $f_1$ achieves its maximum in $r=\ell(1-1/e)$. Then
\begin{equation}
 \frac{F_{\ell}}{({\sf E}X_k)^2}\leq \frac{\exp\left[-\ln n\left(c-\frac{2}{e\ln[1/(1-p)]}\ln\frac{pce}{1-p}+o(1)\right)\right]}{{\sf E}X_k}.
 \label{F_ell_case_2}
\end{equation}
If $c=1$, then $2<\frac{1}{1-p}<3$ and 
$$
\frac{F_{\ell}}{({\sf E}X_k)^2}\leq\exp\left[-\ln n\left(2-\varepsilon-\frac{2}{e\ln[1/(1-p)]}\ln\frac{pe}{1-p}+o(1)\right)\right].
$$ 
Notice that $\ln\frac{pe}{1-p}< 2 \ln\frac{1}{1-p}$. Indeed, $\frac{d}{dx}(2\ln x-\ln(x-1)-1)=\frac{x-2}{x(x-1)}$. Therefore, $2\ln x-\ln(x-1)-1\geq 2\ln 2-1>0$ for $x\geq 2$ (this leads to the above inequality after setting $x=\frac{1}{1-p}$). 

Then,
$$
2-\varepsilon-\frac{2}{e\ln[1/(1-p)]}\ln\frac{pe}{1-p}> 2-\varepsilon-\frac{4}{e}>\frac{1}{2}.
$$

If $c\geq 2$, then $\frac{p}{1-p}<1$. Therefore, $1<\frac{1}{1-p}<2$. Clearly, 
$$
c-\frac{2}{e\ln[1/(1-p)]}\ln\frac{pce}{1-p}\geq
\frac{2}{e\ln[1/(1-p)]}\left(1-\ln\frac{2p}{(1-p)\ln[1/(1-p)]}\right).
$$
Notice that $\ln\frac{2p}{(1-p)\ln[1/(1-p)]}-1<\ln\frac{1}{1-p}$. Indeed, $\frac{d}{dx}(\ln x-\ln(x-1)+\ln\ln x-\ln 2+1)=-\frac{1}{x(x-1)}+\frac{1}{x\ln x}>0$ for $x>1$. Therefore, 
$$
\ln x-\ln(x-1)+\ln\ln x-\ln 2+1=\ln\frac{x\ln(1+(x-1))}{x-1}+\ln\frac{e}{2}>\ln\left(x\left(1-\frac{x-1}{2}\right)\right)+\ln\frac{e}{2}>0
$$ 
for $1<x<2$ (this leads to the above inequality after setting $x=\frac{1}{1-p}$). 

From~(\ref{F_ell_case_2}), we get that
$$
\frac{F_{\ell}}{({\sf E}X_k)^2}\leq \frac{\exp\left[\ln n(\frac{2}{e}+o(1))\right]}{{\sf E}X_k}\leq e^{\ln n\left(\frac{2}{e}-1+\varepsilon+o(1)\right)}.
$$

\item Let $\frac{9}{8}<\frac{pc}{1-p}\leq 2$ and $f_1$ achieves its maximum in $r=\ell/2$. Then
\begin{equation}
 \frac{F_{\ell}}{({\sf E}X_k)^2}\leq \frac{\exp\left[-\ln n\left(c-\frac{1}{\ln[1/(1-p)]}\ln\frac{2pc}{1-p}+o(1)\right)\right]}{{\sf E}X_k}.
\label{F_ell_case_3}
\end{equation}
If $c=1$, then $\frac{17}{8}<\frac{1}{1-p}\leq 3$. Therefore, 
$$
\frac{F_{\ell}}{({\sf E}X_k)^2}\leq\exp\left[-\ln n\left(2-\varepsilon-\frac{\ln[2p/(1-p)]}{\ln[1/(1-p)]}+o(1)\right)\right]=
$$
$$
\exp\left[-\ln n\left(2-\varepsilon-\frac{\ln 2+\ln[1/(1-p)-1]}{\ln[1/(1-p)]}+o(1)\right)\right]<
\exp\left[-\ln n\left(1-\frac{\ln 2}{\ln (17/8)}\right)\right].
$$ 

If $c=2$, then $\frac{25}{16}<\frac{1}{1-p}\leq 2$ and
$$
 \frac{F_{\ell}}{({\sf E}X_k)^2}\leq \frac{\exp\left[-\ln n\left(2-\frac{\ln 4+\ln[1/(1-p)-1]}{\ln[1/(1-p)]}+o(1)\right)\right]}{{\sf E}X_k}.
$$
Since $\frac{d}{dx}\frac{\ln 4+\ln(x-1)}{\ln x}=\frac{x\ln x-(x-1)\ln [4(x-1)]}{x(x-1)\ln^2 x}$; $\frac{d}{dx}(x\ln x-(x-1)\ln[4(x-1)])<0$ on $(25/16,2]$, and $x\ln x-(x-1)\ln [4(x-1)]|_{x=2}=0$, the function $\frac{\ln 4+\ln(x-1)}{\ln x}$ increases on $(25/16,2]$. Therefore, $\frac{\ln 4+\ln[1/(1-p)-1]}{\ln[1/(1-p)]}\leq 2$. Thus, $\frac{F_{\ell}}{({\sf E}X_k)^2}\leq n^{-1+\varepsilon+o(1)}.$

Finally, let $c\geq 3$. Clearly, $p\leq\frac{2}{5}$ and 
$$
c-\frac{1}{\ln[1/(1-p)]}\ln\frac{2pc}{1-p}\geq
\frac{1}{\ln[1/(1-p)]}\left(1-\ln\frac{2p}{(1-p)\ln[1/(1-p)]}\right).
$$
Since $\ln\frac{2p}{(1-p)\ln[1/(1-p)]}$ increases, it is at most $\ln\frac{4}{3\ln(5/3)}<1$. From~(\ref{F_ell_case_3}), we get there exists $\delta>0$ such that, for all large enough $n$,
$\frac{F_{\ell}}{({\sf E}X_k)^2}\leq \frac{\exp\left[-\delta\ln n\right]}{{\sf E}X_k}$.

\item
Finally, assume that $\frac{pc}{1-p}>2$. Then,
$$ 
 \frac{F_{\ell}}{({\sf E}X_k)^2}\leq
 \frac{\exp\left[-\ln n\left(c-\frac{2}{\ln[1/(1-p)]}\ln\frac{pc}{1-p}+o(1)\right)\right]}{{\sf E}X_k}
 $$ 
 
If $c=1$, then $\frac{1}{1-p}>3$. In this case,
$$ 
 \frac{F_{\ell}}{({\sf E}X_k)^2}\leq
 e^{-\ln n\left(2-\varepsilon-\frac{2\ln[1/(1-p)-1]}{\ln[1/(1-p)]}+o(1)\right)}=
 e^{\ln n\left(\varepsilon-\frac{2\ln[1/(1-p)]-2\ln[1/(1-p)-1]}{\ln[1/(1-p)]}+o(1)\right)}<e^{-\varepsilon\ln n(1+o(1))}.
 $$ 
 
If $c=2$, then $\frac{1}{1-p}>2$. Clearly, $x^3>4(x-1)^2+1$ for $x>2$. Therefore,
$$ 
 \frac{F_{\ell}}{({\sf E}X_k)^2}\leq
 e^{-\ln n\left(3-\varepsilon-\frac{\ln(4[1/(1-p)-1]^2)}{\ln[1/(1-p)]}+o(1)\right)}\leq
 e^{\ln n\left(\varepsilon-\frac{\ln\left(1+\frac{1}{4[1/(1-p)-1]^2}\right)}{\ln[1/(1-p)]}+o(1)\right)}<e^{-\varepsilon\ln n(1+o(1))}.
 $$ 
 
Finally, let $c\geq 3$. Then, $\frac{1}{1-p}>\frac{2}{c}+1$. Let us show that, for every $x>\frac{2}{c}+1$, the following inequality holds: $x^{c+1}>c^2(x-1)^2+1$. Indeed,
 $$
  x^{c+1}=((x-1)+1)^{c+1}>1+{c+1\choose 2}(x-1)^2+{c+1\choose 3}(x-1)^3+{c+1\choose 4}(x-1)^4>
 $$
 $$
 1+c(c+1)(x-1)^2\left(\frac{1}{2}+\frac{c-1}{6}\cdot\frac{2}{c}+\frac{(c-1)(c-2)}{24}\cdot\frac{4}{c^2}\right)=
 1+c(c+1)(x-1)^2\left(1-\frac{5}{6c}+\frac{1}{3c^2}\right)=
 $$
 $$
 1+(x-1)^2\left(c^2-\frac{5}{6}c+\frac{1}{3}+c-\frac{5}{6}+\frac{1}{3c}\right)>
 1+(x-1)^2\left(c^2+\frac{1}{6}c-\frac{1}{2}\right)\geq 1+c^2(x-1)^2.
 $$
 Therefore,
 $$ 
 \frac{F_{\ell}}{({\sf E}X_k)^2}\leq
 e^{-\ln n\left(c+1-\varepsilon-\frac{\ln(c^2[1/(1-p)-1]^2)}{\ln[1/(1-p)]}+o(1)\right)}\leq
 e^{\ln n\left(\varepsilon-\frac{\ln\left(1+\frac{1}{c^2[1/(1-p)-1]^2}\right)}{\ln[1/(1-p)]}+o(1)\right)}<e^{-\varepsilon\ln n(1+o(1))}.
 $$ 
  
\end{enumerate}

\section{Maximum induced subgraphs with $t$ edges}
\label{sparse_proof}

Let $\varepsilon>0$ be as small as desired. Set $\varphi(k)=k^2/t(k)$. We know that $\varphi(k)>1/\varepsilon$ for all $k$ large enough. 
	
\subsection{Computing expectation and variance}
	
Let $X_k$ be the number of induced subgraphs in $G(n,p)$ with $k$ vertices and $t(k)$ edges. Then
\begin{equation}
{\sf E}X_k={n\choose k}{{k\choose 2}\choose t}p^t(1-p)^{{k\choose 2}-t},\quad
{\sf E}X_k(X_k-1)={n\choose k}p^{2t}(1-p)^{2{k\choose 2}-2t}\sum_{\ell=0}^{k-1}F_{\ell},
\label{computed_moments}
\end{equation}
where
$$ 
F_{\ell}={k\choose\ell}{n-k\choose k-\ell}(1-p)^{-{\ell\choose 2}}\sum_{j=\max\left\{0,t-{k\choose 2}+{\ell\choose 2}\right\}}^{\min\{t,{\ell\choose 2}\}}{{\ell\choose 2}\choose j}{{k\choose 2}-{\ell\choose 2}\choose t-j}^2\left(\frac{1-p}{p}\right)^j.
$$
	
Everywhere in this section, we assume that $k=\Theta(\ln n)$. Then
$$
{\sf E}X_k=\frac{1}{\sqrt{kt}}\exp\biggl[k\ln n-k\ln k +k+t\ln{k\choose 2}-t\ln t+t+t\ln\frac{p}{1-p}+
$$
$$
\sum_{m=2}^{\infty} (-1)^{m-1}{\frac{t^m}{m({k\choose 2}-t)^{m-1}}}-{k\choose 2}\ln(1/[1-p])+O(1)\biggr]=
$$
$$
=\frac{1}{\sqrt{kt}}\exp\biggl[k\biggl( \ln n+\frac{k}{\varphi(k)}\left[\ln {\frac{{k\choose2}}{t}} +1+\ln \frac{p}{1-p}\right]+
\sum_{m=2}^{\infty}\frac{k(-1)^{m-1}}{m \varphi(k)^m \left(\frac{1}{2}-\frac{1}{\varphi(k)}\right)^{m-1}}-\frac{k}{2}\ln\left(\frac{1}{1-p}\right)+O(\ln k)\biggr) \biggr].
$$
We may assume that $1+\ln\frac{p}{1-p}<\ln \frac{{k\choose2}}{t}< \ln\frac{1}{\varepsilon}$. Since 
$$
\left|\sum_{m=2}^{\infty}(-1)^{m-1}\frac{k}{m \varphi(k)^m (\frac{1}{2}-\frac{1}{\varphi(k)})^{m-1}}\right|<\frac{k}{\varphi(k)}\sum_{m=1}^{\infty}\left[\varphi(k) \left(\frac{1}{2}-\frac{1}{\varphi(k)}\right)\right]^{-m}=\frac{2k}{\varphi(k)(\varphi(k)-4)}<\varepsilon k,
$$
we get
$$
{\sf E}X_k=\frac{1}{\sqrt{kt}}e^{k\left(\ln n -\frac{k}{2}\ln\frac{1}{1-p}+f(k)\right)},
$$
where $|f(k)|<3\varepsilon \ln\frac{1}{\varepsilon}\ln n$.\\
	


	
There exists an $\varepsilon_1\in(0,1/4)$ such that, for all small enough $\varepsilon_0>0$, all $k$ from
\begin{equation}
\left[\left(1-3\varepsilon \ln\frac{1}{\varepsilon}-\varepsilon_0\right)\frac{2}{\ln[1/(1-p)]}\ln n,\left(1+3\varepsilon \ln\frac{1}{\varepsilon}+\varepsilon_0\right)\frac{2}{\ln[1/(1-p)]}\ln n \right]
\label{k_interval}
\end{equation}
and large enough $n$,
$$
\frac{{\sf E}X_{k+1}}{{\sf E}X_k}<
\exp\biggl[\ln n-\ln k-k\ln(1/[1-p])+\left(-\ln t(k) +2 \ln k+1-\ln 2+\ln\frac{p}{1-p}\right) (t(k+1)-t(k))+
$$
$$
\left(\ln {\frac{k^2+k}{k^2-k}}-\ln\frac{t(k+1)}{t(k)}\right)t(k+1)+
\sum_{m=2}^{\infty}(-1)^{m-1}\left(\frac{t^m(k+1)}{m({{k+1}\choose 2}-t(k+1))^{m-1}}- \frac{t^m(k)}{m({{k}\choose 2}-t(k))^{m-1}} \right)+O(1)\biggr]<
$$
$$
e^{-(1-\varepsilon_1)\ln n}.
$$
Indeed,
$$
\left[-\ln t(k) +2 \ln k+1-\ln 2+\ln\frac{p}{1-p}\right] (t(k+1)-t(k))\leq 2R\frac{k\ln\varphi(k)}{\varphi(k)}<2R\varepsilon\ln\frac{1}{\varepsilon}k,
$$
$$
t(k+1)\left[\ln {\frac{k^2+k}{k^2-k}}-\ln\frac{t(k+1)}{t(k)}\right]\leq 
t(k+1)\left[\frac{2}{k-1}+\frac{R}{k}\right]<\frac{(k+1)^2}{\varphi(k+1)}\frac{R+3}{k}<\varepsilon(R+3)(k+3)
$$
and
$$
\frac{t^m(k+1)}{m({{k+1}\choose 2}-t(k+1))^{m-1}}- \frac{t^m(k)}{m({{k}\choose 2}-t(k))^{m-1}}<
$$	
$$
\frac{t^m(k+1)}{m({k\choose 2}-t(k))^{m-1}}\left(1-\frac{2k-2R\varepsilon k}{k^2+k-2t(k+1)}\right)^{m-1}-\frac{t^m(k)}{m({{k}\choose 2}-t(k))^{m-1}}<
$$
$$
\frac{t^m(k+1)-t^m(k)}{m({k\choose 2}-t(k))^{m-1}}<
\frac{R\varepsilon k m(\varepsilon (k+1)^2)^{m-1}}{m({k\choose 2}-\varepsilon k^2)^{m-1}}<2R\varepsilon k\left(\frac{\varepsilon}{1/2-\varepsilon}\right)^{m-1}.
$$

Since, for $k>\frac{2}{\ln(1/[1-p])}\ln n\left[1+3\varepsilon\ln\frac{1}{\varepsilon}+\varepsilon_0\right]$, ${\sf E}X_k\to 0$, and, for $k<\frac{2}{\ln(1/[1-p])}\ln n\left[1-3\varepsilon\ln\frac{1}{\varepsilon}-\varepsilon_0\right]$, ${\sf E}X_k\to \infty$, the minimum $k$ such that ${\sf E}X_k<1$ belongs to~(\ref{k_interval}). If ${\sf E}X_{k-1}>n^{1-2\varepsilon_1}$,  then denote this minimum $k$ by $k_0$. Otherwise, set $k_0=k-1$.

\begin{claim}
Let $C>0$, $k_0+1\leq k<C\ln n$. Then ${\sf E}X_k\to 0$. Moreover, ${\sf E}X_{k-1} >n^{1-2\varepsilon_1}\to\infty$.
\label{cl_1}
\end{claim}

{\it Proof.} If ${\sf E}X_{k_0}<1$, then, by the definition, ${\sf E}X_{k_0-1}>n^{1-2\varepsilon_1}$. Moreover, ${\sf E}X_{k_0+1}\leq n^{\varepsilon_1-1}{\sf E}X_{k_0}<n^{\varepsilon_1-1}\to 0$. Otherwise, ${\sf E}X_{k_0-1}>n^{1-\varepsilon_1}{\sf E}X_{k_0}\geq n^{1-\varepsilon_1}$. Moreover, ${\sf E}X_{k_0+1}\leq n^{\varepsilon_1-1}{\sf E}X_{k_0}\leq n^{-\varepsilon_1}\to 0$.

It remains to notice that if $k>k_0+1$ and belongs to (\ref{k_interval}), then ${\sf E}X_k\leq{\sf E}X_{k_0+1}\to 0$.$\quad\Box$



\subsection{Upper bound}
	
Fix $C>0$ as large as desired. From Markov's inequality and Claim~\ref{cl_1}, it follows that, if $k_0+1\leq k<C\ln n$, then ${\sf P}(X_k\geq 1)\to 0$. Let us prove that the quantification over $k$ can be moved inside the probability.
	
Since $\varepsilon_1<1/4$, if $k_0+3\leq k<C\ln n$, then 
$$
{\sf P}(X_k\geq 1)\leq{\sf E}X_{k}<e^{-2(1-\varepsilon_1)\ln n}{\sf E}X_{k_0+1}<e^{-2(1-\varepsilon_1)\ln n}=o(n^{-3/2}).
$$
	
Now, let $k>C\ln n$. For large enough $n$,
$$
{\sf E}X_k={n\choose k}{{k\choose 2}\choose t}p^t(1-p)^{{k\choose 2}-t}<
n^k t^{-t}e^{t} {k\choose 2}^t p^t (1-p)^{{k\choose 2} -t}=
$$
$$
\exp\left[k\ln n+t\ln{k\choose 2}+t-t\ln t+t\ln\frac{p}{1-p}-{k\choose 2}\ln \frac{1}{1-p}\right]<
$$
$$
\exp\left[\frac{k^2}{C}+t\left(\ln\frac{{k\choose 2}}{t}+\ln\frac{p}{1-p}+1\right)-{k\choose 2}\ln\frac{1}{1-p}\right]<
$$
$$
<\exp\left[\frac{k^2}{C}+2t\ln\frac{1}{\varepsilon}-{k\choose2}\ln\frac{1}{1-p}\right]<\exp\left[k^2\left(\frac{1}{C}+2\varepsilon\ln\frac{1}{\varepsilon}-\frac{1}{2}\ln\frac{1}{1-p}\right)(1+o(1))\right]=o(n^{-3/2}).
$$	
Finally, we get
$$
{\sf P}(\forall k\geq k_0+1\,\, X_k=0)=1-{\sf P}(\exists k\geq k_0+1\,\, X_k\geq 1)\geq 
$$
$$
1-\sum_{k=k_0+1}^n{\sf P}(X_k\geq 1)\geq 1-\sum_{k=k_0+1}^n{\sf E}X_k\geq 1-o(n^{-1/2})\to 1.
$$
	
\subsection{Lower bound}

Here we prove that ${\sf P}(X_{k_0-1}=0)\to 0$.\\
	
Set $k=k_0-1$. Let us estimate ${\sf E}X_k(X_k-1).$
	
Let $A,B$ be $k$-vertex subsets of $\{1,\ldots,n\}$ having at most 1 common vertex. For a $k$-vertex set $W\subset\{1,\ldots,n\}$, let $I[W]=1$, if $G(n,p)|_W$ has exactly $t$ edges, and let $I[W]=0$ otherwise. Clearly, ${\sf E}I[A]I[B]={\sf E}I[A]{\sf E}I[B]$. Therefore, $F_0+F_1-({\sf E}X_k)^2\leq 0$. Then
\begin{equation}
{\sf D}X_k={\sf E}X_k(X_k-1)+{\sf E}X_k-({\sf E}X_k)^2\leq {n\choose k} p^{2t}(1-p)^{2{k\choose 2}-2t}\sum_{\ell=2}^{k-1}F_{\ell}+{\sf E}X_k.
\label{DXk}
\end{equation}
Let us estimate $F_{\ell}$.
	
\subsubsection{Small $\ell$}
	
Let 
$$
\delta\in\left(\frac{1/(1-p)}{\ln[1/(1-p)]}\cdot\frac{4\varepsilon}{1-2\varepsilon},\frac{\ln[1/(1-p)]}{12p/(1-p)+\ln[1/(1-p)]}\right).
$$
Consider
$$
\ell\leq\ell^*:=\left\lfloor \frac{2-\delta}{\ln[1/(1-p)]}\ln n\right\rfloor.
$$

By the definition of $F_{\ell}$,
\begin{equation}
{n\choose k}p^{2t}(1-p)^{2{k\choose 2}-2t}\frac{F_{\ell}}{({\sf E}X_k)^2}=\frac{{k\choose\ell}{n-k\choose k-\ell}}{{n\choose k}}(1-p)^{-{\ell\choose 2}}\sum_{j=0}^{\min\{t,{\ell\choose 2}\}}
\frac{{{\ell\choose 2}\choose j}{{k\choose 2}-{\ell\choose 2}\choose t-j}^2}{{{k\choose 2}\choose t}^2}\left(\frac{1-p}{p}\right)^j.
\label{F_ell_new_expression}
\end{equation}

If $j<t$, then
$$
\frac{{{k\choose 2}-{\ell\choose 2}\choose t-j}^2}{{{k\choose 2}\choose t}^2}\leq
\frac{{{k\choose 2}\choose t-j}^2}{{{k\choose 2}\choose t}^2}\leq
\frac{e^4}{2\pi^2}\times\frac{{k\choose 2}-t}{{k\choose 2}-(t-j)}\times\frac{t}{t-j}\times
$$
$$
\left(1+\frac{j}{t-j}\right)^{2(t-j)}t^{2j}
\left(1-\frac{j}{{k\choose2}-(t-j)}\right)^{2\left({k\choose2}-(t-j)\right)}\left({k\choose2}-t\right)^{-2j}<
\frac{e^4}{2\pi^2}t\left(\frac{t}{{k\choose 2}-t}\right)^{2j}.
$$
If $j=t$, then 
$$
\frac{{{k\choose 2}-{\ell\choose 2}\choose t-j}^2}{{{k\choose 2}\choose t}^2}=\frac{1}{{{k\choose 2}\choose t}^2}\leq \left [ \frac{e^2({k\choose 2}-t)^{{k\choose 2}-t+1/2}t^{t+1/2}}{{\sqrt{2\pi}} {k\choose 2}^{{k\choose 2}-t+1/2}{k\choose 2}^t}\right] ^2\leq\frac{e^4}{2\pi}t\left(\frac{t}{{k\choose 2}}\right)^{2t}.
$$
Summing up, from~(\ref{F_ell_new_expression}), we get
$$
{n\choose k}p^{2t}(1-p)^{2{k\choose 2}-2t}\frac{F_{\ell}}{({\sf E}X_k)^2}\leq G_{\ell}\sum_{j=0}^{\ell\choose 2} {{\ell\choose 2}\choose j} \frac{e^4}{2\pi} {t}\left(\frac{t}{{k\choose 2}-t}\right)^{2j}\left(\frac{1-p}{p}\right)^j=
$$
\begin{equation} 
\frac{e^4}{2\pi}G_{\ell}{t}\left(1+\frac{1-p}{p}\left(\frac{t}{{k\choose 2}-t}\right)^{2}\right)^{{\ell\choose 2}}=\frac{e^4}{2\pi}\widetilde{G}_{\ell}{t},
\label{est0}
\end{equation}
where 
\begin{equation}
G_{\ell}=\frac{{k\choose\ell}{n-k\choose k-\ell}}{{n\choose k}}(1-p)^{-{\ell\choose 2}},\quad
\tilde{G}_{\ell}=G_{\ell}(A(t))^{\ell\choose 2},\quad A(t)=1+\frac{1-p}{p}\left(\frac{t}{{k\choose 2}-t}\right)^{2}.
\label{G_def}
\end{equation}

Clearly, $\frac{\tilde{G}_{\ell+1}}{\tilde{G}_{\ell}}=\frac{(k-\ell)^2}{(\ell+1)(n-2k+\ell+1)}(1-p)^{-\ell} A(t)^{\ell}$. Therefore,
\begin{equation}
\frac{\partial}{\partial\ell}\ln\frac{\tilde{G}_{\ell+1}}{\tilde{G}_{\ell}}=-\frac{2}{k-\ell}-\frac{1}{\ell+1}-\frac{1}{n-2k+\ell+1}+\ln\left(\frac{1}{1-p}\right)+\ln A(t).
\label{tilde_G_derivative}
\end{equation}
	
Since $\frac{1}{1-p}>1$ is a constant, for certain non-negative $\delta_1(n)=O(1)$, $\delta_2(n)=O(1)$, the right side of~(\ref{tilde_G_derivative}) is positive when $\ell\in(\delta_1,k-\delta_2)$ and negative when $\ell\in(0,\delta_1)$ and $\ell\in(k-\delta_2,k)$. Since $\tilde G_3/\tilde G_2=O(k^2/n)$ (less than 1 for large enough  $n$), $\tilde G_{k-1}/\tilde G_{k-2}>n^{1-o(1)}$ (bigger than 1 for large enough $n$), on $[2,k-2]$, there exists a unique $\ell_0$ (not necessarily integer) such that 

$\tilde G_{\ell_0+1}/\tilde G_{\ell_0}=1$, 

if $\ell<\ell_0$, then $\frac{\tilde{G}_{\ell+1}}{\tilde{G}_{\ell}}<1$, 

if $\ell>\ell_0$, then $\frac{\tilde G_{\ell+1}}{\tilde{G}_{\ell}}>1$. 

So, $\tilde{G}_{\ell}$ as function of integer argument $\ell$ first decreases, and then increases. Therefore, for $\ell\in\{2,3,\ldots,\ell_*\}$,
$$
\tilde{G}_{\ell}\leq\max\{\tilde G_2,\tilde G_{\ell_*}\}\leq\max\left\{\frac{k^2{n\choose k-2}}{(1-p){n\choose k}}A(t),G_{\ell_*}A(t)^{{\ell_*\choose 2}}\right\}.
$$
The first value equals $O\left(\frac{k^4}{n^2}\right)$. Let us estimate the second value:
$$
G_{\ell_*}A(t)^{{\ell_*\choose 2}}\leq\frac{k^{\ell_*}(n-k)^{k-\ell_*}k^k}{(n-k)^k}\left[\frac{1}{1-p}e^{\frac{1-p}{p}\left(\frac{t}{{k\choose 2}-t}\right)^2}\right]^{(\ell_*)^2/2}<k^{2k}\left[n^{-\delta/2+\frac{1/(1-p)}{\ln[1/(1-p)]}\cdot\frac{2\varepsilon}{1-2\varepsilon}(1-\delta/2)+o(1)}\right]^{\ell_*}<
$$
$$ 
\exp\left[-\left(\delta/2-\frac{1/(1-p)}{\ln[1/(1-p)]}\cdot\frac{2\varepsilon}{1-2\varepsilon}+o(1)\right)\frac{2-\delta}{\ln\frac{1}{1-p}}\ln^2 n\right].
$$

From~(\ref{est0}), we get
\begin{equation}
{n\choose k}p^{2t}(1-p)^{2{k\choose 2}-2t}\frac{F_{\ell}}{({\sf E}X_k)^2}=O(k^6/n^2).
\label{small_ell_bound}
\end{equation}

\subsubsection{Large $\ell$}
	
Let $\ell> \ell_*$. Denote ${\ell\choose 2 }=L$, ${k\choose 2}=K$. 

Clearly,
\begin{equation}
{n\choose k}p^{2t}(1-p)^{2K-2t}\frac{F_{\ell}}{({\sf E}X_k)^2}= 
\frac{{k\choose\ell}{n-k\choose k-\ell}(1-p)^{K-L}}{{\sf E}X_k}\sum_{j=\max\left\{0,t-K+L\right\}}^{t}H_{\ell,j},
\label{est}
\end{equation}
where 
\begin{equation}
H_{\ell,j}=\frac{{L\choose j}{K-L\choose t-j}^2}{{K\choose t}}\left(\frac{1-p}{p}\right)^{j-t}.
\label{H_definition}
\end{equation}
If $j=t$, then
$$
H_{\ell,j}=\frac{{L\choose t}}{{K\choose t}}<1.
$$
If $j=t-K+L$, then
$$
H_{\ell,j}=\frac{{L\choose t-K+L}}{{K\choose t}}\left(\frac{1-p}{p}\right)^{-K+L}<\left(\frac{tp}{L(1-p)}\right)^{K-L}<1
$$
as well. Therefore, in both cases, we get
\begin{equation}
\frac{{k\choose\ell}{n-k\choose k-\ell}(1-p)^{K-L}}{{\sf E}X_k} H_{\ell,j}<\frac{{k\choose\ell}{n-k\choose k-\ell}(1-p)^{K-L}}{{\sf E}X_k}<
\frac{1}{{\sf E}X_k}e^{(k-\ell)\left[\ln\frac{ke}{k-\ell}+\ln n-\frac{k+\ell+1}{2}\ln\frac{1}{1-p}\right]}<
\frac{1}{{\sf E}X_k}e^{(k-\ell)\left[\ln\frac{ke}{\ell}-\frac{1}{2}\ln n\right]}
\label{boarder_j}
\end{equation}
since $\ell>k/2$.

If $j=0$, then
$$
H_{\ell,j}=\frac{{K-L\choose t}^2}{{K\choose t}}\left(\frac{1-p}{p}\right)^{-t}=
\frac{(K-L-t+1)^2 p}{Lt^2(1-p)}\frac{L{K-L\choose t-1}^2}{{K\choose t}}\left(\frac{1-p}{p}\right)^{1-t}=
$$
\begin{equation}
\frac{(K-L-t+1)^2 p}{Lt^2(1-p)}H_{\ell,1}<\frac{k^2p}{1-p}H_{\ell,1}.
\label{j=0}
\end{equation}


Finally, let $\max\{t-K+L,0\}+1\leq j\leq t-1$. Clearly,
$$
\frac{{L\choose j}{{K-L}\choose {t-j}}}{{K\choose t}}=
O\left( \sqrt{t}\frac{L^L(K-L)^{K-L}t^t(K-t)^{K-t}}{j^j (L-j)^{L-j}(t-j)^{t-j}(K-L-t+j)^{K-L-t+j}K^K}\right)=
$$
$$
O\left( \sqrt{t}\frac{\left(\frac{L}{j}\right)^j\left(1+\frac{j}{L-j}\right)^{L-j}\left(\frac{K-L}{t-j}\right)^{t-j}\left(1+\frac{t-j}{K-L-t+j}\right)^{K-L-t+j}}{{{{(K/t)}^{t}\left(1+\frac{t}{K-t}\right)^{K-t}}}}\right).
$$
Let us prove that, for positive $a,b,A,B$ such that $a<A, b<B$, the following is true: $\left(1+\frac{a}{A}\right)^A\left( 1+\frac{b}{B}\right)^B<\left(1+\frac{a+b}{A+B}\right)^{A+B}$. Consider the fraction 
$$
\left.\left(1+\frac{a}{A}\right)^A\left( 1+\frac{b}{B}\right)^B\right/\left(1+\frac{a+b}{A+B}\right)^{A+B}=
$$
$$
=\left(1+\frac{aB-Ab}{A(A+B+a+b)}\right)^A\left(1-\frac{aB-Ab}{B(A+B+a+b)}\right)^B<
e^{\frac{aB-Ab}{A+B+a+b}}e^{-\frac{aB-Ab}{A+B+a+b}}=1.
$$
Then $\left(1+\frac{j}{L-j}\right)^{L-j}\left(1+\frac{t-j}{K-L-t-j}\right)^{K-L-t+j}<\left(1+\frac{t}{K-t}\right)^{K-t}$, 
and so
$$
H_{\ell,j}=O\left(\sqrt{t}\frac{\left( \frac{L}{j}\right)^j\left(\frac{K-L}{t-j}\right)^{t-j}}{(K/t)^{t}}\right){{K-L}\choose{t-j}}\left(\frac{1-p}{p}\right)^{j-t}=
$$
$$
=O\left(\sqrt{t}e^{j\ln L-j\ln j+2(t-j)\ln\left(K-L\right)-2(t-j)\ln(t-j)+(t-j)-t\ln K+t\ln t+(j-t)\ln (\frac{1-p}{p})}\right)
$$
$$
\leq \mathrm{exp}\biggl(2j\ln k-j\ln 2-j\ln j+2(t-j)\ln(k-\ell)+2(t-j)\ln k-2(t-j)\ln(t-j)-
$$
\begin{equation}
2t\ln k+t\ln 2+t\ln t+
(j-t)\left[\ln \left(\frac{1-p}{p}\right)-1\right]+O(\ln n)\biggr)=
e^{h(j)+t\ln 2+t\ln t+O(\ln n)},
\label{H_upper}
\end{equation}
where 
$$
h(j)=-j\ln 2-j\ln j+2(t-j)\ln(k-\ell)-2(t-j)\ln(t-j)+(j-t)\left[\ln \left(\frac{1-p}{p}\right)-1\right].
$$
Compute the derivative
$$
\frac{\partial h}{\partial j}=
\ln \frac{1-p}{p}-\ln 2-\ln j-2\ln(k-\ell)+2\ln(t-j) =\ln\frac{(t-j)^2}{j}-\ln\frac{2p(k-\ell)^2}{1-p}.
$$
It equals 0 if and only if
\begin{equation}
\frac{(t-j)^2}{j}=\frac{2p(k-\ell)^2}{1-p}.
\label{roots}
\end{equation} 
If the last equality is true, then
\begin{equation}
h(j)=t-j-t\ln 2-t\ln j+(t-j)\ln\frac{2p(k-\ell)^2}{1-p}-(t-j)\ln\frac{(t-j)^2}{j}=
t-j-t\ln 2-t\ln j
\label{in_j_1}
\end{equation}
and so
\begin{equation}
h(j)+t\ln 2+t\ln t=t-j-t\ln (j/t)=t-j-t\ln\left(1-\frac{t-j}{t}\right).
\label{in_maxima}
\end{equation}
Let us estimate both roots $j_1,j_2$ of the equation~(\ref{roots}). Clearly, $h$ increases when $j<j_1$, decreases when $j_1<j<j_2$ and again increases when $j>j_2$. The roots equal
$$
j_{1,2}=t+\frac{p}{1-p}(k-\ell)^2\left(1\pm\sqrt{1+\frac{2t(1-p)}{p(k-\ell)^2}}\right).
$$
Then, $j_2>t$ and $t>j_1>0$ ($j_1< t$ since $\sqrt{1+\frac{2t(1-p)}{p(k-\ell)^2}}> 1$; $j_1> 0$ since the product of the roots equals $t^2>0$). Therefore, the maximum of $h(j)$ on $[0,t-1]$ is at most $h(j_1)$.

Now, let us finish an upper bound for the left side of~(\ref{in_maxima}). 

First, assume that $t\geq\frac{p}{2(1-p)}(k-\ell)^2$. Then $\frac{t-j_1}{t}<\frac{\sqrt{2}-1}{2}$ since the function $x(\sqrt{1+2t/x}-1)$ increases in $x$. Therefore, $t\ln(1-\frac{t-j_1}{t})>\frac{p\ln[(3-\sqrt{2})/2]}{2(1-p)}(k-\ell)^2$.

Second, assume that $t<\frac{p}{2(1-p)}(k-\ell)^2$. Then
$$
 \frac{t-j_1}{t}<\frac{p(k-\ell)^2}{t(1-p)}\left(1+\frac{t(1-p)}{p(k-\ell)^2}-\frac{1}{2}\left(\frac{t(1-p)}{p(k-\ell)^2}\right)^2+\frac{1}{2}\left(\frac{t(1-p)}{p(k-\ell)^2}\right)^3-1\right)=
$$
$$
 1-\frac{t(1-p)}{2p(k-\ell)^2}+\frac{1}{2}\left(\frac{t(1-p)}{2p(k-\ell)^2}\right)^2.
$$
Therefore,
$$
 t\ln\left(1-\frac{t-j_1}{t}\right)>t\ln\left(\frac{t(1-p)}{2p(k-\ell)^2}-\frac{1}{2}\left(\frac{t(1-p)}{2p(k-\ell)^2}\right)^2\right)>t\ln\frac{t(1-p)}{4p(k-\ell)^2}>-\frac{4p(k-\ell)^2}{1-p}.
$$

Finally, we get that
$$
h(j)+t\ln 2+t\ln t\leq 2(t-j_1)-t\ln\left(1-\frac{t-j_1}{t}\right)\leq 
2\frac{p}{1-p}(k-\ell)^2\sqrt{1+\frac{2t(1-p)}{p(k-\ell)^2}}+\frac{4p}{1-p}(k-\ell)^2\leq 
$$
$$ 
\frac{6p}{1-p}(k-\ell)^2+2\frac{\sqrt{2p}}{\sqrt{1-p}}\sqrt{t}(k-\ell).
$$
	
From this and (\ref{DXk}),~(\ref{G_def}),~(\ref{small_ell_bound}),~(\ref{est}),~(\ref{boarder_j}),~(\ref{j=0}),~(\ref{H_upper}), we get that
$$
{\sf P}(X_k=0)\leq\frac{{\sf D}X_k}{({\sf E}X_k)^2}\leq\frac{1+o(1)}{{\sf E}X_k}+O\left(\frac{k^6}{n^2}\right)+{n\choose k} p^{2t}(1-p)^{2{k\choose 2}-2t}\sum_{\ell=\ell^*+1}^{k-1}\frac{F_{\ell}}{({\sf E}X_k)^2}\leq
$$
$$
\frac{(1-p)^{k\choose 2}{n\choose k}}{{\sf E}X_k}\sum_{\ell=\ell^*+1}^{k-1}t\frac{k^2p}{1-p}G_{\ell}e^{2\frac{\sqrt{2p}}{\sqrt{1-p}}\sqrt{t}(k-\ell)+\frac{6p}{1-p}(k-\ell)^2}+o(1).
$$
Let $\hat G_{\ell}=G_{\ell}e^{2\frac{\sqrt{2p}}{\sqrt{1-p}}\sqrt{t}(k-\ell)+\frac{6p}{1-p}(k-\ell)^2}$. Obviously, 
$$
\frac{\partial}{\partial\ell}\ln\frac{\hat G_{\ell+1}}{\hat G_{\ell}}=\frac{\partial}{\partial\ell}\ln\frac{G_{\ell+1}}{G_{\ell}}+\frac{12p}{1-p}=-\frac{2}{k-\ell}-\frac{1}{\ell+1}-\frac{1}{n-2k+\ell+1}+\ln\left(\frac{1}{1-p}\right)+\frac{12p}{1-p}.
$$
Therefore, in the same way as for $\tilde G_{\ell}$, on $[2,k-2]$ there exists a unique $\ell_0$ (not necessarily integer) such that 

$\hat G_{\ell_0+1}/ \hat G_{\ell_0}=1$, 

if $\ell<\ell_0$, then $\frac{\hat G_{\ell+1}}{\hat G_{\ell}}<1$,

if $\ell>\ell_0$, then $\frac{\hat G_{\ell+1}}{\hat G_{\ell}}>1$. 

Let us show that $\ell_0<\ell^*$. 


Clearly,
$$
1=\hat G_{\ell_0+1} / \hat G_{\ell_0}=\frac{(k-\ell_0)^2}{(\ell_0+1)(n-2k+\ell_0+1)}(1-p)^{-\ell_0}e^{-2\frac{\sqrt{2p}}{\sqrt{1-p}}\sqrt{t}-\frac{6p}{1-p}(2k-2\ell_0-1)}.
$$
Then, $\ell_0\ln\frac{1}{1-p}-2\frac{\sqrt{2p}}{\sqrt{1-p}}\sqrt{t}-\frac{6p}{1-p}(2k-2\ell_0)=\ln n(1+o(1))$. Therefore, 
$$
\ell_0=\frac{\ln n}{\ln\frac{1}{1-p}}\left(2-\frac{\ln\frac{1}{1-p}-4\frac{\sqrt{2p}}{\sqrt{1-p}}\sqrt{\varepsilon}}{\ln\frac{1}{1-p}+\frac{12p}{1-p}}+o(1)\right)<\ell^*.
$$

Therefore, $\hat G_{\ell}$ increases if $\ell>\ell^*$.\\
	
Finally, we get
$$
{\sf P}(X_k=0) \leq \frac{tk^3p(1-p)^{{k\choose 2}-1}{n\choose k}}{{\sf E}X_k}G_{k-1}e^{2\frac{\sqrt{2p}}{\sqrt{1-p}}\sqrt{t}+\frac{6p}{1-p}}+o(1)=\frac{tk^4(n-k)p(1-p)^{k-2}}{{\sf E}X_k}e^{2\frac{\sqrt{2p}}{\sqrt{1-p}}\sqrt{t}+\frac{6p}{1-p}}+o(1)\leq 
$$
$$
n^{-1+o(1)}e^{2\frac{\sqrt{2p}}{\sqrt{1-p}}\sqrt{t}}+o(1)=o(1).
$$

\section*{Acknowledgments} 
The research is supported by the grant 16-11-10014 of Russian Science Foundation.


\begin{thebibliography}{9}

\bibitem{BE_independent}  Bollob\'{a}s B. and Erd\H{o}s P., {\em Cliques in random graphs}, Math. Proc. Camb. Phil. Soc. {\bf 80} (1976), 419--427.

\bibitem{Grimmett_independent} G.R. Grimmett, C.J.H. McDiarmid, {\it On colouring random graphs}, Math. Proc. Camb. Phil. Soc. {\bf 77} (1975) 313--324.

\bibitem{M_independent_0} D. Matula, {\it The employee party problem}, Not. Amer. Math. Soc., 19(2): A--382, 1972.

\bibitem{M_independent} Matula D., {\em The largest clique size in a random graph}, Tech. Rep. Dept. Comp. Sci., Southern Methodist University, Dallas, Texas, 1976.



\bibitem{Bollobas} Bollob\'{a}s B.{\em Random Graphs}, 2nd Edition, Cambridge University Press, 2001.

\bibitem{Janson} Janson S., \L uczak T. and Ruci\'{n}ski A., {\em Random Graphs}, New York, Wiley, 2000.

\bibitem{Survey} Raigorodskii A.M., Zhukovskii M.E., {\em Random graphs: models and asymptotic characteristics}, Russian Mathematical Surveys {\bf 70}:1 (2015)~33--81.

\bibitem{ZhukRodi} Rodionov I.V., Zhukovskii M.E, {\it On the Distribution of the Maximum k-Degrees of the Binomial Random Graph}, Doklady Mathematics {\bf 98}:3 (2018)~619--621.


\bibitem{Many_Induced_Graphs} Dutta K. and Subramanian C.R., {\em On Induced Paths, Holes and Trees in Random Graphs}, Proc. ANALCO 2018, 168--177.

\bibitem{Subgraphs_bounded_deg} Fountoulakis N., Kang R.J., McDiarmid C., {\em Largest sparse subgraphs of random graphs}, European Journal of Combinatorics {\bf 35} (2014), 232--244.

\bibitem{BaloghZhuk} J. Balogh, M. Zhukovskii, {\it On the sizes of large subgraphs of the binomial random graph}, 2019, https://arxiv.org/pdf/1904.05307.pdf.

\bibitem{DerevyankoKiselev} Derevyanko N.M., Kiselev S.G., {\it Independence Numbers of Random Subgraphs of Some Distance Graph}, Problems of Information Transmission, {\bf 53}:4 (2017) 307--318.


\bibitem{RaiStability} Raigorodskii A.M., {\it On the stability of the independence number of a random subgraph}, Doklady Mathematics, {\bf 96}:3 (2017) 628--630.







\end{thebibliography}
\end{document}